\documentclass[titlepage,11pt]{article}
\usepackage{latexsym}
\usepackage{amsfonts}
\usepackage{amsmath}
\newcommand\blackslug{\hbox{\hskip 1pt \vrule width 4pt height 8pt depth 1.5pt
        \hskip 1pt}}
\newcommand\bbox{\hfill \quad \blackslug \bigbreak}

\def\dd{\hbox{-}}
\def\cc{\hbox{-}\cdots\hbox{-}}
\def\ll{,\ldots,}

\title{Induced subgraphs of graphs with large chromatic number.
\\VIII. Long odd holes}
\author{Maria Chudnovsky\thanks{Supported by NSF grant DMS-1550991 and US Army Research Office Grant W911NF-16-1-0404.}\\
Princeton University, Princeton, NJ 08544, USA
\\
\\
Alex Scott\\
Mathematical Institute, University of Oxford, Oxford OX2 6GG, UK
\\
\\
Paul Seymour\thanks{Supported by ONR grant N00014-14-1-0084 and NSF
grant DMS-1265563.}\\
Princeton University, Princeton, NJ 08544, USA
\\
\\
Sophie Spirkl\thanks{Current address: Rutgers University, New Brunswick, NJ 08901, USA}\\
Princeton University, Princeton, NJ 08544, USA}

\date{December 12, 2016; revised \today}

\newtheorem{thm}{}[section]

\newcommand{\Proof}{\noindent{\bf Proof.}\ \ }

\begin{document}
\maketitle
\begin{abstract}
We prove a conjecture of Andr\'as Gy\'arf\'as, that for all $\kappa,\ell$,
every graph with clique number at most $\kappa$ and sufficiently large chromatic number
has an odd hole of length at least $\ell$.
\end{abstract}

\section{Introduction}
All graphs in this paper are finite and have no loops or parallel edges. We denote the chromatic number of a graph $G$
by $\chi(G)$, and its clique number (the cardinality of its largest clique) by $\omega(G)$. A {\em hole} in $G$ means an induced subgraph
which is a cycle of length at least four, and an {\em odd} hole is one with odd length. In~\cite{gyarfas}, Andr\'as Gy\'arf\'as
proposed three conjectures about the lengths of holes in graphs with large chromatic number and 
bounded clique number, the following:

\begin{thm}\label{oddholes} 
For all $\kappa\ge 0$ there exists $c\ge 0$ such that for every graph $G$, if $\omega(G)\le \kappa$
and $\chi(G)>c$ then $G$ has an odd hole.
\end{thm}
\begin{thm}\label{longholes} For all $\kappa,\ell\ge 0$ there exists $c$ such that
for every graph $G$, if $\omega(G)\le \kappa$
and $\chi(G)>c$, then $G$ has a hole of length at least $\ell$.
\end{thm}
\begin{thm}\label{mainthm} For all $\kappa,\ell\ge 0$ there exists $c$ such that
for every graph $G$, if $\omega(G)\le \kappa$
and $\chi(G)>c$ then $G$ has an odd hole of length at least $\ell$.
\end{thm}
The third evidently contains the other two. The first was proved in~\cite{oddholes}, and the second in~\cite{longholes}, but until
now the third has remained open. In this paper we prove the third.

Since this paper was submitted for publication, two of us have proved an even stronger 
result~\cite{holeparity}:
\begin{thm}\label{holeparity}
For all integers $p,q,\kappa\ge 0$ (with $q>0$), there exists $c$ such that
for every graph $G$, if $\omega(G)\le \kappa$
and $\chi(G)>c$, then $G$ has a hole of length $p$ modulo $q$.
\end{thm}

There have been two other partial results approaching \ref{mainthm}, 
in~\cite{pentagonal, holeseq}.
The second is the stronger, namely:
\begin{thm}\label{holeseq}
For all $\ell\ge 0$,
there exists $c\ge 0$ such that for every graph $G$, if $\omega(G)\le 2$
and $\chi(G)>c$ then $G$ has holes of $\ell$ consecutive lengths (and in particular has an odd hole of length at least $\ell$).
\end{thm}
One might expect that an analogous result holds for all graphs $G$ with $\omega(G)\le \kappa$, for all fixed $\kappa$, as was conjectured in~\cite{holeseq}, 
but this remains open.
We can prove that there are two holes both of length at least $\ell$ and with consecutive lengths, by a refinement of the 
methods of this paper, but we omit the details.

The proof of \ref{mainthm} will be by induction on $\kappa$, with $\ell$ fixed, so we may assume that $\kappa\ge 2$
and the result holds for all smaller $\kappa$. In particular
there exists $\tau$ such that for every graph $G$, if $\omega(G)\le \kappa-1$
and $G$ has no odd hole of length at least $\ell$ then $\chi(G)<\tau$. 
We might as well assume that $\ell$ is odd and $\ell\ge 5$.
Let us say a graph $G$ is a {\em $(\kappa,\ell,\tau)$-candidate} if $\omega(G)\le \kappa$, and $G$ has 
no odd hole of length  at least $\ell$, and every induced subgraph of $G$ with clique number less than $\kappa$ is $\tau$-colourable. 
We will show that for all $\kappa,\ell,\tau$ there exists $c$ such that every $(\kappa,\ell,\tau)$-candidate has chromatic number
at most $c$.

If $H$ is a subgraph of $G$, and not necessarily induced, then adjacency in $H$ may be different from in $G$, and we speak of
being $G$-adjacent, a $G$-neighbour, and so on, to indicate in which graph we are 
using adjacency when there may be confusion.
If $X\subseteq V(G)$, the subgraph of $G$ induced on $X$ is denoted by $G[X]$,
and we often write $\chi(X)$ for $\chi(G[X])$. The {\em distance} or {\em $G$-distance} between two vertices $u,v$
of $G$ is the length of a shortest path between $u,v$, or $\infty$ if there is no such path.
If $v\in V(G)$ and $\rho\ge 0$ is an integer, $N_G^{\rho}(v)$ or $N^{\rho}(v)$ denotes the set of all vertices $u$ with distance
exactly
$\rho$  from $v$, and $N_G^{\rho}[v]$ or $N^{\rho}[v]$ denotes the set of all $v$ with distance at most $\rho$ from $v$.
If $G$ is a nonnull graph  and $\rho\ge 1$,
we define $\chi^{\rho}(G)$ to be the maximum of $\chi(N^{\rho}[v])$ taken over all vertices $v$ of $G$.
(For the null graph $G$ we define $\chi^{\rho}(G)=0$.)

As in several other papers of this series, the proof of \ref{mainthm} examines whether there is an induced subgraph
of large chromatic number such that every ball of small radius in it has bounded chromatic number. The proof 
breaks into three steps, as follows. We first show that for every $(\kappa,\ell,\tau)$-candidate $G$, $\chi(G)$ is bounded by a function of $\chi^3(G)$. 
Second, a similar argument shows the same for $\chi^2(G)$; that is, there is a function $\phi$ such that
$\chi(G)\le \phi(\chi^2(G))$ for every $(\kappa,\ell,\tau)$-candidate $G$. 
For the third step,
we apply
a result of~\cite{chandeliers}, which implies that in these circumstances there is an upper bound on the chromatic number of all 
$(\kappa,\ell,\tau)$-candidates.

\section{Gradings}

Let $G$ be a graph. We say a {\em grading} of $G$ is a sequence $(W_1\ll W_n)$ of subsets of $V(G)$, pairwise
disjoint and with union $V(G)$. If $w\ge 0$ is such that $\chi(G[W_i])\le w$ for $1\le i\le n$ we say the grading is {\em $w$-colourable}. We say that $u\in V(G)$ is {\em earlier} than $v\in V(G)$
(with respect to some grading $(W_1\ll W_n)$) if $u\in W_i$ and $v\in W_j$ where $i<j$.
We need some lemmas about gradings. The first is:

\begin{thm}\label{bigtouch}
Let $w,c\ge 0$ and let $(W_1\ll W_n)$ be a $w$-colourable grading of a graph $G$ with $\chi(G)>w+2c$.
Then there exist subsets $X,Y$ of $V(G)$ with the following properties:
\begin{itemize}
\item $G[X], G[Y]$ are both connected;
\item every vertex in $Y$ is earlier than every vertex in $X$;
\item some vertex in $X$ has a neighbour in $Y$; and
\item $\chi(X), \chi(Y)>c$.
\end{itemize}
\end{thm}
\Proof
Let us say $v\in V(G)$ is {\em right-active} if 
there exists $X\subseteq V(G)$ such that
\begin{itemize}
\item $G[X]$ is connected;
\item $v$ is earlier than every vertex in $X$;
\item $v$ has a neighbour in $X$; and
\item $\chi(X)>c$.
\end{itemize}
Let $A$ be the set of all vertices of $G$ that are not right-active, and $B=V(G)\setminus A$.
Since $\chi(G)>w+2c$, either $\chi(A)>w+c$ or $\chi(B)>c$. Suppose that $\chi(A)>w+c$, and let
$C$ be a component of $G[A]$ with maximum chromatic number. Choose $i$
with $1\le i\le n$ minimum such that $W_i\cap V(C)\ne \emptyset$. Since $\chi(W_i)\le w$, it follows that $\chi(C\setminus W_i)>c$.
Let $C'$ be a component of $C\setminus W_i$ with maximum chromatic number; consequently $\chi(C')>c$. Choose $v\in W_i\cap V(C)$
with a neighbour in $C'$ (this exists since $C$ is connected). But then $v$ is right-active with $X=C'$, a contradiction.

This proves that $\chi(A)\le w+c$, and so $\chi(B)>c$. Let $C$ be a component of $G[B]$ with maximum chromatic number. Choose $i$
with $1\le i\le n$ maximum such that $W_i\cap V(C)\ne \emptyset$, and choose $v\in W_i\cap V(C)$. Since $v$ is right-active, there
exists $X$ as in the definition of ``right-active'' above; and then setting $Y=V(C)$ satisfies the theorem. This proves 
\ref{bigtouch}.~\bbox

Next, we need:

\begin{thm}\label{greentouch}
Let $w,c\ge 0$ and let $(W_1\ll W_n)$ be a $w$-colourable grading of a graph $G$. Let $H$ be a subgraph of $G$ (not necessarily induced)
with $\chi(H)>w+2(c+\chi^1(G))$. Then there is an edge $uv$ of $H$, and a subset $X$ of $V(G)$,
such that
\begin{itemize}
\item $G[X]$ is connected;
\item $u,v$ are both earlier than every vertex in $X$;
\item exactly one of $u,v$ has a $G$-neighbour in $X$; and
\item $\chi(G[X])>c$.
\end{itemize}
\end{thm}
\Proof
We claim:
\\
\\
(1) {\em There exist $X,Y\subseteq V(G)$ with the following properties:
\begin{itemize}
\item $G[X]$ and $H[Y]$ are connected;
\item every vertex of $Y$ is earlier than every vertex of $X$;
\item some vertex of $Y$ has no $G$-neighbour in $X$, and some vertex of $Y$ has a $G$-neighbour in $X$; and
\item $\chi(G[X])>c$.
\end{itemize}
}
\noindent By \ref{bigtouch} applied to $H$, there exist $X,Y\subseteq V(H)$ such that 
\begin{itemize}
\item $H[X], H[Y]$ are both connected;
\item every vertex in $Y$ is earlier than every vertex in $X$;
\item some vertex in $X$ has an $H$-neighbour in $Y$; and
\item $\chi(H[X]), \chi(H[Y])>c+\chi^1(G)$.
\end{itemize}
If some vertex of $Y$ has no $G$-neighbour in $X$, then (1) holds, so 
we may 
assume that every vertex of $Y$ has a $G$-neighbour in $X$. Choose $y\in Y$ and let $N$ be the set of vertices in $X$ that are
$G$-adjacent to $y$. Let $C$ be a component of $G[X\setminus N]$ with maximum chromatic number. Since $\chi(G[N])\le \chi^1(G)$,
it follows that $\chi(C)>c$. If some vertex of $Y$ has a $G$-neighbour in $V(C)$ then (1) holds with $X=C$, so we assume not. Choose $x\in N$
with a $G$-neighbour in $V(C)$. Then $G[V(C)\cup \{x\}]$ is connected, and some vertex in $Y$ has a $G$-neighbour in it, namely $y$.
We may therefore assume that every vertex of $Y$ has a $G$-neighbour in $V(C)\cup \{x\}$, and hence is $G$-adjacent to $x$, since
no vertex of $Y$ has a $G$-neighbour in $V(C)$. But this is impossible since $\chi(G[Y])\ge \chi(H[Y]) >\chi^1(G)$. This proves (1).

\bigskip
Let $X,Y$ be as in (1). Since some vertex of $Y$ has a $G$-neighbour in $X$ and some vertex of $Y$ has no $G$-neighbour in $X$,
and $H[Y]$ is connected, it follows that there is an edge $uv$ of $H$ with $u,v\in Y$ such that exactly one of $u,v$
has a $G$-neighbour in $X$; and so the theorem holds. This proves \ref{greentouch}.~\bbox

We also need the following:

\begin{thm}\label{getpath}
Let $G$ be a graph, let $k\ge 0$, let $C\subseteq V(G)$, and let $x_0\in V(G)\setminus C$,
such that
$G[C]$ is connected,
$x_0$ has a neighbour in $C$, and
$\chi(C)>k \chi^1(G)$.
Then there is an induced path $x_0 \cc x_{k}$ of $G$ where $x_1\ll x_{k}\in C$, and a subset $C'$ of $C$,
with the following properties:
\begin{itemize}
\item $x_0\ll x_{k}\notin C'$;
\item $G[C']$ is connected;
\item $x_{k}$ has a neighbour in $C'$, and $x_0\ll x_{k-1}$ have no neighbours in $C'$; and
\item $\chi(C')\ge \chi(C)-k\chi^1(G)$.
\end{itemize}
\end{thm}
\Proof We proceed by induction on $k$; the result holds if $k=0$, so we assume that $k>0$ and the result holds for $k -1$.
Consequently there is an induced path $x_0\cc x_{k-1}$ of $G$ where $x_1\ll x_{k-1}\in C$, and a subset $C''$ of $C$,
such that
\begin{itemize}
\item $x_0\ll x_{k-1}\notin C''$;
\item $G[C'']$ is connected;
\item $x_{k-1}$ has a neighbour in $C''$, and $x_0\ll x_{k-2}$ have no neighbours in $C''$; and
\item $\chi(C'')\ge \chi(C)-(k-1)\chi^1(G)$.
\end{itemize}
Let $N$ be the set of neighbours of $x_{k-1}$, and let $C'$ be the vertex set of a
component of $G[C''\setminus N]$, chosen with $\chi (C')$ maximum (there is such a component
since $\chi(C'')>\chi^1(G)\ge \chi(N)$).
Let $x_{k}$ be a neighbour of $x_{k-1}$ with a neighbour in $C'$. Then $x_0\cc x_{k}$ and $C'$ satisfy the theorem.
This proves \ref{getpath}.~\bbox

If $G$ is a graph and $B,C\subseteq V(G)$, we say that $B$ {\em covers} $C$ if $B\cap C=\emptyset$ and every vertex in $C$ 
has a neighbour in $B$.
Let $G$ be a graph, and let $B,C\subseteq V(G)$, where $B$ covers $C$. Let $B=\{b_1\ll b_m\}$.
For $1\le i<j\le m$ we say that $b_i$ is {\em earlier} than $b_j$
(with respect to the enumeration $(b_1\ll b_m)$). For $v\in C$, let $i\in \{1\ll m\}$ be minimum such that $b_i,v$ are adjacent;
we call $b_i$ the {\em earliest parent} of $v$. An edge $uv$ of $G[C]$ is said to be {\em square} (with respect to the enumeration 
$(b_1\ll b_m)$) if
the earliest parent of $u$ is nonadjacent to $v$, and the earliest parent of $v$ is nonadjacent to $u$.
Let $B=\{b_1\ll b_m\}$, and let $(W_1\ll W_n)$ be a grading of $G[C]$. We say the enumeration $(b_1\ll b_m)$  of $B$ and the 
grading $(W_1\ll W_n)$ are {\em compatible} if for all $u,v\in C$ with $u$ earlier than $v$, the earliest parent of $u$ is earlier than
the earliest parent of $v$.

\begin{thm}\label{getgreenedge}
Let $G$ be a graph, and let $B,C\subseteq V(G)$, where $B$ covers $C$. Let every induced subgraph $J$ of $G$ with $\omega(J)<\omega(G)$
have chromatic number at most $\tau$.
Let the enumeration $(b_1\ll b_m)$ of $B$ and the grading $(W_1\ll W_n)$ of $G[C]$ be compatible. Let $H$ be the subgraph of $G$
with vertex set $C$ and edge set the set of all square edges. Let $(W_1\ll W_n)$ be $w$-colourable; then
$\chi(G[C])\le w\tau\chi(H)$.
\end{thm}
\Proof
Let $X\subseteq C$ be a stable set of $H$. We claim that $\chi(G[X])\le w\tau$. Suppose not. Since $(W_1\ll W_n)$ is $w$-colourable,
there is a partition $(A_1\ll A_w)$ of $C$ such that $W_i\cap A_j$ is stable in $G$ for $1\le i\le n$ and $1\le j\le w$.
Consequently there exists $j$ such that $\chi(G[X\cap A_j])>\tau$. From the choice of $\tau$, it follows that 
$\omega(G[X\cap A_j])=\omega(G)$. Let $Y\subseteq X\cap A_j$ be a clique with $|Y|=\omega(G)$. 
Choose $i\in \{1\ll n\}$ maximum such that
$W_i\cap Y\ne \emptyset$, and let $y\in W_i\cap Y$. Let $b$ be the earliest parent of $y$. Since $|Y|=\omega(G)$, there exists $y'\in Y$
nonadjacent to $b$. Since $y,y'\in Y\subseteq A_j$, and $y,y'$ are adjacent, it follows that $y'\notin W_i$, and so $y'$ is earlier 
than $y$, from the choice of $y$. Since the enumeration and the grading are compatible, the earliest parent $b'$ of $y'$ is 
earlier than the earliest parent of $y$, and in particular is nonadjacent to $y$; 
and so $yy'$ is a square edge, contradicting that $y,y'\in X$. This proves that $\chi(G[X])\le w\tau$.
Since $C$ can be partitioned into $\chi(H)$ sets that are stable in $H$, it follows that $\chi(G[C])\le w\tau\chi(H)$.
This proves \ref{getgreenedge}.~\bbox

Combining these lemmas yields the main result of this section:
\begin{thm}\label{mainlemma}
Let $G$ be a graph, and let $B,C\subseteq V(G)$, where $B$ covers $C$. Let every induced subgraph $J$ of $G$ with $\omega(J)<\omega(G)$
have chromatic number at most $\tau$. Let $(W_1\ll W_n)$ be a $w$-colourable grading of $G[C]$, and 
let the enumeration $(b_1\ll b_m)$ of $B$ be compatible with this grading. Let $\ell,\rho\ge 1$
be integers. Let 
$$\chi(C)>w\tau (4(\ell+2)\chi^{\rho}(G)+w).$$
Then there is an induced path $p_1\cc p_k$ of $G[C]$, such that 
\begin{itemize}
\item $k\ge \ell$;
\item $p_1p_2$ is a square edge;
\item $p_1,p_2$ are both earlier than all of $p_3\ll p_k$; and
\item let $b,b'$ be the earliest parent of $p_1,p_2$ respectively; then for each $v\in \{b,b', p_1\ll p_{\ell}\}$,
the $G$-distance between $p_k$ and $v$ is at least $\rho+1$.
\end{itemize}
\end{thm}
\Proof 
Let $c=\ell\chi^1(G)+(\ell+2)\chi^{\rho}(G)$. Thus 
$$\chi(C)>w\tau (2(c+\chi^1(G))+w).$$
Let $H$ be the subgraph of $G[C]$ with edge set the square edges; then by \ref{getgreenedge}, 
$$\chi(H)> 2(c+\chi^1(G))+w.$$
By \ref{greentouch}, applied to $G[C]$ and $H$, there is a square edge $uv$, and a subset $X$ of $C$,
such that
\begin{itemize}
\item $G[X]$ is connected;
\item $u,v$ are both earlier than every vertex in $X$;
\item exactly one of $u,v$ is adjacent in $G$ to a member of $X$; and
\item $\chi(G[X])>c$.
\end{itemize}
Let $\{p_1,p_2\}=\{u,v\}$ where $p_2$ has a neighbour in $X$. Since $\chi(X)>c\ge \ell\chi^1(G)$, \ref{getpath} implies that
there is an induced path $p_1 \cc p_{\ell}$ of $G$ where $p_3\ll p_{\ell}\in X$, and a subset $C'$ of $X$,
with the following properties:
\begin{itemize}
\item $p_1\ll p_{\ell}\notin C'$;
\item $G[C']$ is connected;
\item $p_{\ell}$ has a neighbour in $C'$, and $p_1\ll p_{\ell-1}$ have no neighbours in $C'$; and
\item $\chi(C')\ge \chi(X)-\ell\chi^1(G)> (\ell+2)\chi^{\rho}(G)$.
\end{itemize}
Let $b,b'$ be the earliest parent of $p_1,p_2$ respectively, and let $Z=\{b,b',p_1\ll p_{\ell}\}$. Thus $|Z|=\ell+2$.
Since $\chi(C')>(\ell+2)\chi^{\rho}(G)$, there exists $x\in C'$ such that the $G$-distance between $x$ and each member of $Z$
is at least $\rho+1$. Let $p_{\ell}\cc p_k$ be an induced path of $G[C'\cup \{p_{\ell}\}]$ between $p_{\ell}$ and $x$; then 
$p_1\cc p_k$ satisfies the theorem. This proves \ref{mainlemma}.~\bbox

\section{Little bounding balls}

In this section we carry out the first two steps of the proof, showing that the chromatic number of every
$(\kappa,\ell,\tau)$-candidate $G$ can be bounded in terms of $\chi^3(G)$, and then can be bounded in terms of $\chi^2(G)$. 
More precisely we first prove the following.
\begin{thm}\label{3balls}
Let $G$ be a $(\kappa,\ell,\tau)$-candidate; then $\chi(G)\le 24(2\ell+5)\tau\chi^3(G)^2$.
\end{thm}
\Proof
We may assume that $\chi(G)>4\chi^3(G)$.
Since some component of $G$ has the same chromatic number, we may assume that $G$ is connected.
Choose some vertex, and for $i\ge 0$ let $L_i$ be the set of vertices with $G$-distance $i$ from the vertex. There exists $s$ such that 
$\chi(L_{s+1})\ge \chi(G)/2$. Since $\chi(L_0\cup \cdots\cup L_3)\le \chi^3(G)$, it follows that $s\ge 3$.
Let $V_0$ be the vertex set of a component of $G[L_{s+1}]$ with maximum chromatic number; and choose $z\in L_s$ with a neighbour in $V_0$. 
For $i\ge 0$
let $M_i$ be the set of vertices in $V_0\cup \{z\}$ with $G[V_0\cup \{z\}]$-distance $i$ from $z$. 
Choose $t$ such that $\chi(M_{t+1})\ge \chi(V_0)/2$.
Again, since $\chi(V_0)/2> \chi^3(G)$, it follows that $t\ge 3$. 

The set of vertices in $M_{t+1}$ with $G$-distance at most three from $z$ has chromatic number at most $\chi^3(G)$;
so there is a set $C\subseteq M_{t+1}$ with $\chi(C)\ge \chi(M_{t+1})-\chi^3(G)$ such that 
every vertex in $C$ has $G$-distance at least four from $z$. Let $B$ be the set of vertices in $L_{s}$ with a neighbour in $C$,
and let $D$ be the set of vertices in $M_t$ with a neighbour in $C$. Thus every vertex in $B\cup D$ has $G$-distance at least three
from $z$.

Let $V_1=M_0\cup\cdots\cup M_{t-1}$. Thus $G[V_1]$ is connected, and there are no edges between $V_1$ and $C$.
Let $B_0$ be the set of vertices in $B$ with no neighbour in $V_1$. Let $B_1$ be the set of vertices $v$ in $B$ with a neighbour in $V_1$
such that the $G[V_1\cup \{v\}]$-distance between $z,v$ is odd, and let $B_2$ be the set where this distance is even. Every vertex in $C$
has a neighbour in at least one of $B_0,B_1,B_2$; let $C_i$ be the set of vertices in $C$ with a neighbour in $B_i$ for $i = 0,1,2$.
\\
\\
(1) {\em $\chi(C_0)\le (4\ell+9)\tau\chi^3(G)^2$.}
\\
\\
Take an enumeration $(d_1\ll d_n)$ of $D$, and for $1\le i\le n$ let $W_i$ be the set of vertices $v\in C_0$
such that $v$ is adjacent to $d_i$ and nonadjacent to $d_1\ll d_{i-1}$.
Then $(W_1\ll W_n)$ is a grading of $C_0$, and is $\chi^3(G)$-colourable (indeed, $\tau$-colourable, but we need the $\chi^3(G)$ bound),
and the enumeration $(d_1\ll d_n)$ is compatible with it.
Suppose that $\chi(C_0)>(4\ell+9)\tau\chi^3(G)^2$.
By \ref{mainlemma} with $\rho=3$ and $w=\chi^3(G)$, there is an induced path $p_1\cc p_k$ of $G[C_0]$, such that
\begin{itemize}
\item $k\ge \ell$;
\item $p_1p_2$ is a square edge;
\item $p_1,p_2$ are both earlier than all of $p_3\ll p_k$; and
\item let $d,d'$ be the earliest parent of $p_1,p_2$ respectively; then for each $v\in \{d,d', p_1\ll p_{\ell}\}$,
the $G$-distance between $p_k$ and $v$ is at least $4$.
\end{itemize}
Since $d,d'\in D$, there are induced paths $Q,Q'$ of $G[V_1\cup D]$ between $d,z$ and between $d',z$ respectively, both of length $t$. 
Let $y\in B_0$ be adjacent to $p_k$. 
There is an induced path $R$ between $z,y$ with interior in $L_0\cup\cdots\cup L_{s-1}$, since $z,y\in L_s$. Now $y$ is nonadjacent to
$d,d'$ since the $G$-distance between $p_k$ and $d,d'$ is at least four. It follows that $Q\cup R$ is an induced path between $d,y$, and $Q'\cup R$
is an induced path between $d',y$, of the same length. Also $y$ is nonadjacent to $p_1,p_2\ll p_{\ell}$ since $p_k$ has $G$-distance four
from all these vertices. Choose $j\le k$ minimum such that $p_j$ is adjacent to $y$. Now $d,d'$ are both nonadjacent to $p_3\ll p_{j}$ since $p_1, p_2$
are both earlier than $p_3\ll p_{j}$; and since $p_1$ is nonadjacent to $d'$ and $p_2$ is nonadjacent to $d$ 
(because $p_1p_2$ is a square edge)
it follows that $d\dd p_1\dd p_2\cc p_j\dd y$ and $d'\dd p_2\cc p_j\dd y$ are both induced paths, joining $d,y$ and $d',y$ respectively. So the union of
$d\dd p_1\dd p_2\cc p_j\dd y$ with $Q\cup R$ is a hole, and the union of $d'\dd p_2\cc p_j\dd y$ with $Q'\cup R$ is a hole; and these holes differ in length by one.
Since they both have length more than $\ell$, this is impossible. This proves (1).
\\
\\
(2) {\em For $h = 1,2$, $\chi(C_h)\le (4\ell+9)\tau\chi^3(G)^2.$}
\\
\\
Enumerate the vertices of $V_1$ in increasing order of $G[V_1]$-distance from $z$, breaking ties arbitrarily; that is,
take an enumeration $(a_1\ll a_n)$ of $V_1$ where for $0\le i<j\le n$ the $G[V_1]$-distance between $a_i,z$ in $G[V_1]$
is at most that between $a_j,z$.
For each $v\in C_h$, choose $i\in \{1\ll n\}$ minimum such that some vertex in $B_h$ is adjacent both to $v$ and to $a_i$.
(Such a value of $i$ exists from the definition of $B_h,C_h$.) We call $a_i$ the {\em earliest grandparent} of $v$.
For $1\le i\le n$, let $W_i$ be the set of vertices in $C_h$ with earliest grandparent $a_i$. Thus every vertex in $W_i$
has $G$-distance two from $a_i$, and so $\chi(W_i)\le \chi^3(G)$; and it follows that $(W_1\ll W_n)$ is a $\chi^3(G)$-colourable
grading of $G[C_h]$ (indeed, it is $\chi^2(G)$-colourable, 
since $W_i\subseteq N^2[a_i]$). 
Take an enumeration $(b_1\ll b_m)$ of $B_h$ such that vertices with earlier neighbours in $V_1$ come first;
that is, such that for $1\le i <j\le m$, there is a neighbour $a_h\in V_1$
of $b_i$ such that $h\le h'$ for every neighbour $a_{h'}$ of $b_j$ in $V_1$. It follows that $(b_1\ll b_m)$ and $(W_1\ll W_n)$
are compatible.

Suppose that $\chi(C_h)>(4\ell+9)\tau\chi^3(G)^2$. By \ref{mainlemma} with $\rho=3$ and $w=\chi^3(G)$,
there is an induced path $p_1\cc p_k$ of $G[C_h]$, such that
\begin{itemize}
\item $k\ge \ell$;
\item $p_1p_2$ is a square edge;
\item $p_1,p_2$ are both earlier than all of $p_3\ll p_k$; and
\item let $b,b'$ be the earliest parent of $p_1,p_2$ respectively; then for each $v\in \{b,b', p_1\ll p_{\ell}\}$,
the $G$-distance between $p_k$ and $v$ is at least $4$.
\end{itemize}
Let $b''$ be the earliest parent of $p_k$, and let $a,a',a''$ be the earliest grandparents of $p_1,p_2, p_k$ respectively.
It follows that $ab,a'b', a''b''$ are edges. Since both $a,a'$ occur in the enumeration $(a_1\ll a_n)$ before $a''$, and $a''$ is the
earliest grandparent of $p_k$, there are induced paths $Q,Q'$ of $G[V_1]$ between $a,z$ and between $a', z$ respectively, with lengths of the same parity as
$h+1$
(because $p_1,p_2\in C_h$), such that
$b''$ has no neighbours in $Q,Q'$. There is an induced path $R$ between 
$z,b''$ with interior in $L_0\cup\cdots\cup L_{s-1}$, using only
two vertices of $L_{s-1}$ (neighbours of $z,b''$ respectively). Now since $p,p'\in C$, they both have 
$G$-distance at least four from $z$; and so $b,b'$ both have
$G$-distance at least three from $z$. Moreover $b,b'$ both have $G$-distance at least four from $p_k$ and hence at least three
from $b''$. 
Consequently they both have $G$-distance at least two from the vertices of $R$ in $L_{s-1}$ (it is to arrange this
that we need to control the chromatic number of balls of radius three); and so $b,b'$ both have no neighbours in $R$.
It follows that $Q\cup R$ is an induced path between $b,b''$, and $Q'\cup R$
is an induced path between $b',b''$, and they have lengths of the same parity. Now $b''$ is nonadjacent to $p_1\ll p_{\ell}$ since the $G$-distance between $p_k$
and $p_1\ll p_{\ell}$ is at least four. Choose $j$ minimum such that $p_j$ is adjacent to $b''$. Also, $b,b'$ are nonadjacent to $p_3\ll p_k$
since $p_1,p_2$ are earlier than $p_3\ll p_k$; and since $p_1p_2$ is square, the union of $Q\cup R$ with $b\dd p_1\dd p_2\cc p_j\dd b''$
and the union of $Q'\cup R$ with $b'\dd p_2\cc p_j\dd b''$ are holes of opposite parity, both of length more than $\ell$, which is impossible. This
proves (2).

\bigskip
From (1) and (2), we deduce that 
$\chi(C)\le 3(4\ell+9)\tau\chi^3(G)^2$. But $\chi(C)\ge \chi(M_{t+1})-\chi^3(G)$, so $\chi(M_{t+1})\le 3(4\ell+9)\tau\chi^3(G)^2+\chi^3(G)$.
Since $\chi(M_{t+1})\ge \chi(V_0)/2$, and $\chi(V_0)\ge \chi(G)/2$, we deduce that 
$$\chi(G)\le 12(4\ell+9)\tau\chi^3(G)^2+ 4\chi^3(G)\le 24(2\ell+5)\tau\chi^3(G)^2.$$
This proves \ref{3balls}.~\bbox

Next we bound $\chi(G)$ in terms of $\chi^2(G)$, as follows.

\begin{thm}\label{2balls}
Let $G$ be a $(\kappa,\ell,\tau)$-candidate; then $\chi(G)\le 96(2\ell+5)^3\tau^3\chi^{2}(G)^8$.
\end{thm}
\Proof
Let $z$ be a vertex such that $\chi(N^3[z])=\chi^3(G)$, and let $L_i$ be the set of vertices with $G$-distance $i$ from $z$, for $0\le i\le 3$.
Fix a $\chi^2(G)$-colouring of $G[L_0\cup L_1\cup L_2]$, and for each $v\in L_3$ choose a path of length three from $z$ to $v$, and 
let $\alpha(v), \beta(v)$ be the colours of its second and third vertex respectively. Choose colours $\alpha,\beta$, and let $C$ 
be the set of $v\in L_3$ such that $\alpha(v)=\alpha$ and $\beta(v)=\beta$.
Let $B$ be the set of vertices in $L_2$ with colour $\beta$ and with a neighbour in $L_1$
with colour $\alpha$. Consequently $B$ covers $C$, and any two vertices in $B$ are joined by an induced path of even length with interior in $L_0\cup L_1$.

Let $(b_1\ll b_n)$ be some enumeration of $B$, and for $1\le i\le n$ let $W_i$ be the set of $v\in C$ such that $v$ is adjacent to $b_i$ but not to $b_1\ll b_{i-1}$. 
Thus $(W_1\ll W_n)$ is a $\tau$-colourable grading of $G[C]$ compatible with $(b_1\ll b_n)$. Suppose that $\chi(C)> (4\ell+9)\tau\chi^{2}(G)^2$. 
Then by \ref{mainlemma} with $\rho=2$ and $w=\chi^2(G)$
there is an induced path $p_1\cc p_k$ of $G[C]$, such that
\begin{itemize}
\item $k\ge \ell$;
\item $p_1p_2$ is a square edge;
\item $p_1,p_2$ are both earlier than all of $p_3\ll p_k$; and
\item let $b,b'$ be the earliest parent of $p_1,p_2$ respectively; then for each $v\in \{b,b', p_1\ll p_{\ell}\}$,
the $G$-distance between $p_k$ and $v$ is at least $3$.
\end{itemize}
It follows that $b,b'$ are nonadjacent to $p_3\ll p_k$; let $b''\in B$ be adjacent to $p_k$, and choose $j$ minimum such that $p_j$ is adjacent to $b''$.
Since $p_k$ has $G$-distance at least three from each of $b,b', p_1\ll p_{\ell}$, it follows that $b''$ is nonadjacent to all these vertices. Choose $j$
minimum such that $p_j$ is adjacent to $b''$; then $b\dd p_1\dd p_2\dd p_j\dd b''$ and $b'\dd p_2\cc p_j\dd b''$ are induced paths both of length more than $\ell$
and with lengths of opposite parity. Let $Q$ be one of them with odd length. There is an induced path of even length joining the ends of $Q$
with interior in $L_0\cup L_1$; and its union with $Q$ is an odd hole of length more than $\ell$, which is impossible.

This proves that $\chi(C)\le (4\ell+9)\tau\chi^{2}(G)^2$. Since this holds for every choice of $\alpha,\beta$, and there are only $\chi^2(G)^2$ such choices,
it follows that 
$\chi(L_3)\le (4\ell+9)\tau\chi^{2}(G)^4.$
But $\chi^3(G)\le \chi(L_3)+\chi^2(G)$, so
$\chi^3(G)\le 2(2\ell+5)\tau\chi^{2}(G)^4$.
From \ref{3balls}, it follows that
$\chi(G)\le 96(2\ell+5)^3\tau^3\chi^{2}(G)^8$. This proves \ref{2balls}.~\bbox

\section{Multicovers}

In this section we combine \ref{2balls} with a result of~\cite{chandeliers} to deduce \ref{mainthm}, and for this we need 
some definitions.
If $X,Y$ are disjoint subsets of the vertex set of a graph $G$, we say
\begin{itemize}
\item $X$ is {\em complete} to $Y$ if every vertex in $X$
is adjacent to every vertex in $Y$; and
\item $X$ is {\em anticomplete} to $Y$ if every vertex in $X$
nonadjacent to every vertex in $Y$.
\end{itemize}
(If $X=\{v\}$ we say $v$ is complete to $Y$ instead of $\{v\}$, and so on.)

Let $x\in V(G)$, let $N$ be some set of neighbours of $x$, and let $C\subseteq V(G)$
be disjoint from $N\cup \{x\}$, such that $x$ is anticomplete to $C$ and $N$ covers $C$. In this situation
we call $(x,N)$ a {\em cover} of $C$ in $G$. For $C,X\subseteq V(G)$, a {\em multicover of $C$} in $G$
is a family $(N_x:x\in X)$
such that
\begin{itemize}
\item $X$ is stable;
\item for each $x\in X, (x,N_x)$ is a cover of $C$;
\item for all distinct $x,x'\in X$, $x'$ is anticomplete to $N_x$ (and in particular all the sets
$\{x\}\cup N_x$ are pairwise disjoint).
\end{itemize}
Its {\em length} is $|X|$; and the
multicover $(N_x:x\in X)$ is {\em stable} if each of the sets $N_x\;(x\in X)$ is stable.

\begin{thm}\label{getstablemult}
For all $m,c,\kappa,\tau\ge 0$, suppose that $\omega(G)\le \kappa$, and $\chi(H)\le \tau$ for every induced
subgraph $H$ of $G$ with $\omega(H)<\kappa$. If 
$G$ admits a multicover $(N_x:x\in X)$ with length $m$  of
a set $C'\subseteq V(G)$
with $\chi(C')>c\tau^m$,
then $G$ admits a stable multicover contained in $(N_x:x\in X)$ with length $m$, of some subset $C\subseteq C'$
with $\chi(C)>c$.
\end{thm}
\Proof
Let
$(N_x':x\in X)$ be a multicover in $G$ of length $m$, of a set $C'$ with $\chi(C')> c\tau^m$.
For each $x\in X$, since $G[N_x']$ has clique number less than $\kappa$, 
this subgraph is 
$\tau$-colourable; choose some such colouring, 
with colours $1\ll \tau$ (for each $x$).
For each $v\in C'$, let $f_v:X\rightarrow \{1\ll \tau\}$ such that for each $x\in X$, some neighbour of $v$ in $N_x'$
has colour $f_v(x)$. There are only $\tau^m$ possibilities for $f_v$, so there is a function $f:X\rightarrow \{1\ll \tau\}$
and a subset $C\subseteq C'$ with $\chi(C)\ge \chi(C')\tau^{-m}> c$, such that $f_v=f$ for all $v\in C$. For each $x\in X$,
let $N_x$ be the set of vertices in $N_x'$ with colour $f(x)$; then $(N_x:x\in X)$ is a stable multicover of $C$.
This proves \ref{getstablemult}.~\bbox

Let $(N_x:x\in X)$ be a multicover of $C$ in $G$. It is said to be {\em $k$-crested} if there are vertices $a_1\ll a_k$ and vertices $a_{ix}(1\le i\le k, x\in X)$ of $G$,
all distinct, with the following properties:
\begin{itemize}
\item $a_1\ll a_k$ and the vertices $a_{ix}(1\le i\le k, x\in X)$ do not belong to $C\cup X\cup \bigcup_{x\in X}N_x$;
\item for $1\le i\le k$  and each $x\in X$, $a_{ix}$ is adjacent to $x$, and there are no other edges between the sets 
$\{a_1\ll a_k\}\cup \{a_{ix}:1\le i\le k, x\in X\}$ and $C\cup X\cup \bigcup_{x\in X}N_x$;
\item for $1\le i\le k$  and each $x\in X$, $a_{ix}$ is adjacent to $a_i$, and there are no other edges between $\{a_1\ll a_k\}$ and $\{a_{ix}:1\le i\le k, x\in
X\}$
\item $a_1\ll a_k$ are pairwise nonadjacent;
\item for all distinct $i,j\in \{1\ll k\}$ and all $x,y\in X$, $a_{ix}$ is nonadjacent to $a_{jy}$.
\end{itemize}
Note that $a_{ix}$ may be adjacent to $a_{iy}$. We say the multicover is {\em stably $k$-crested} if for $1\le i\le k$ and all distinct $x,y\in X$, $a_{ix},
a_{iy}$ are nonadjacent.
Theorem 2.1 of~\cite{longholes} (setting $j$ of the theorem to be $\kappa$) implies:

\begin{thm}\label{gettick}
For all $m,c,\kappa,\tau\ge 0$ there exist $m',c'\ge 0$ with the following property.
Let $G$ be a graph with $\omega(G)\le \kappa$, such that $\chi(H)\le \tau$ for every induced subgraph $H$ of $G$ with $\omega(H)<\kappa$.
Let $(N_x':x\in X')$ be a stable multicover in $G$ of some set $C'$, such that $|X'|\ge m'$ and $\chi(C')> c'$.
Then there exist $X\subseteq X'$ with $|X|\ge m$, and $C\subseteq C'$ with $\chi(C)> c$, and a stable multicover $(N_x:x\in X)$
of $C$ contained in $(N_x':x\in X')$ that is 1-crested.
\end{thm}
Ramsey's theorem applied to the vertices $\{a_{1x}:x\in X\}$ together with \ref{gettick} yields (under the same hypotheses) that such a stable multicover
exists which is stably 1-crested; and combining this with \ref{getstablemult} yields:
\begin{thm}\label{getbettertick}
For all $m,c,\kappa,\tau\ge 0$ there exist $m',c'\ge 0$ with the following property.
Let $G$ be a graph with $\omega(G)\le \kappa$, such that $\chi(H)\le \tau$ for every induced subgraph $H$ of $G$ with $\omega(H)<\kappa$.
Let $(N_x':x\in X')$ be a multicover in $G$ of some set $C'$, such that $|X'|\ge m'$ and $\chi(C')> c'$.
Then there exist $X\subseteq X'$ with $|X|\ge m$, and $C\subseteq C'$ with $\chi(C)> c$, and a stable multicover $(N_x:x\in X)$
of $C$ contained in $(N_x':x\in X')$ that is stably 1-crested.
\end{thm}

Repeated application of this yields:
\begin{thm}\label{getbigtick}
For all $m,c,k,\kappa,\tau\ge 0$ there exist $m',c'\ge 0$ with the following property.
Let $G$ be a graph with $\omega(G)\le \kappa$, such that $\chi(H)\le \tau$ for every induced subgraph $H$ of $G$ with $\omega(H)<\kappa$.
Let $(N_x':x\in X')$ be a  multicover in $G$ of some set $C'$, such that $|X'|\ge m'$ and $\chi(C')> c'$.
Then there exist $X\subseteq X'$ with $|X|\ge m$, and $C\subseteq C'$ with $\chi(C)> c$, and a stable multicover $(N_x:x\in X)$
of $C$ contained in $(N_x':x\in X')$ that is stably $k$-crested.
\end{thm}

We deduce:
\begin{thm}\label{cresttohole}
For all $\ell,\kappa,\tau\ge 0$, there exist $m',c'\ge 0$ such that if $G$ satisfies:
\begin{itemize}
\item $\omega(G)\le \kappa$;
\item $\chi(H)\le \tau$ for every induced subgraph $H$ of $G$ with $\omega(H)<\kappa$; and
\item $G$ admits a multicover of length at least $m'$ of a set $C\subseteq V(G)$ with $\chi(C)>c'$;
\end{itemize}
then $G$ has an odd hole of length at least $\ell$.
\end{thm}
\Proof Let $m',c'$ satisfy \ref{getbigtick}, with $k,m$ both replaced by $\ell$ and with $c=\tau$.
Let $G$ be as in the theorem; so $G$ admits a stably
$\ell$-crested stable multicover $(N_x:x\in X)$ of a set $C$,
where $|X|=\ell$ and $\chi(C)>c$. Let $a_1\ll a_{\ell}$ and the vertices $a_{ix}(1\le i\le \ell, x\in X)$ be as in the definition of stable $k$-crested.
Choose distinct $x_1,x_2\in X$. Let $\ell_1\in \{\ell+1,\ell+3\}$ be a multiple of four. There is an induced path $P$ between $x_1,x_2$ of length $\ell_1$
such that $V(P)\subseteq \{a_1\ll a_{\ell}\}\cup X\cup \{a_{ix}:1\le i\le \ell, x\in X\}$. Since $\chi(C)>\tau$, 
there is a $\kappa$-clique $Y\subseteq C$.
Choose $b_1\in N_{x_1}\cup N_{x_2}$ with as many neighbours in $Y$ as possible. By exchanging $x_1,x_2$
if necessary, we may assume that $b_1\in N_{x_1}$. Now $b_1$ is not complete to $Y$
since $\omega(G)\le \kappa$; so there exists $y_2\in Y$ nonadjacent to $b_1$. Choose $b_2\in N_{x_2}$ adjacent to $y_2$. From the choice of $b_1$, there exists
$y_1\in Y$ adjacent to $b_1$ and not to $b_2$. If $b_1,b_2$ are adjacent let $Q$ be the path $x_1\dd b_1\dd b_2\dd x_2$, and otherwise let $Q$
be the path $x_1\dd b_1\dd y_1\dd y_2\dd b_2\dd x_2$. Thus $Q$ is induced, and has length three or five, and the 
union of $P$ and $Q$ is an odd hole of length more than $\ell$. This proves \ref{cresttohole}.~\bbox

Let $\mathbb{N}$ denote the set of nonnegative integers, let $\phi:\mathbb{N}\rightarrow \mathbb{N}$ be a nondecreasing function, and
let $h\ge 1$ be an integer. We say a graph $G$ is {\em $(2,\phi)$-controlled} if 
$\chi(H)\le \phi(\chi^2(H))$ for every induced subgraph $H$ of $G$.
Consequently \ref{2balls} implies that:
\begin{thm}\label{1-control}
Every $(\kappa,\ell,\tau)$-candidate $G$ is $(2,\phi)$-controlled where $\phi$ is the function
$$\phi(x) = 96(2\ell+5)^3\tau^3x^8.$$
\end{thm}

We need the following, a consequence of theorem 9.7 of~\cite{chandeliers}. That involves ``trees of lamps'', but we 
do not need to define 
those here; all we need is that a cycle of length $\ell$ is a tree of lamps. (Note that what we call a
``multicover'' here is called a ``strongly-independent 2-multicover'' in that paper, and indexed in a slightly different way.)
\begin{thm}\label{findchand5}
Let $m,\kappa,c',\ell\ge 0$, and let $\phi:\mathbb{N}\rightarrow \mathbb{N}$ 
be non-decreasing. Then there exists $c$ with the following property. Let $G$ be a graph
such that
\begin{itemize}
\item $\omega(G)\le \kappa$;
\item $G$ is $(2,\phi)$-controlled;
\item $G$ does not admit a multicover of length $m$ of a set with chromatic number more than $c'$; and
\item $G$ has no hole of length $\ell$.
\end{itemize}
Then $\chi(G)\le c$.
\end{thm}

\bigskip

Finally, we can prove \ref{mainthm}, which we restate:
\begin{thm}\label{mainthmagain}
For all $\kappa,\ell\ge 0$ there exists $c$ such that
for every graph $G$, if $\omega(G)\le \kappa$
and $\chi(G)>c$ then $G$ has an odd hole of length at least $\ell$.
\end{thm}
\Proof
As we saw in section 1, it suffices to prove an upper bound on the chromatic number of all $(\kappa,\ell,\tau)$-candidates,
for all $\tau\ge 0$. Let $m',c'$ satisfy \ref{cresttohole}, and let 
$\phi$ be as in \ref{1-control}. Let $c$ satisfy \ref{findchand5}, with $m$ replaced by $m'$.

Now let $G$ be a $(\kappa,\ell,\tau)$-candidate.
If $G$ admits a multicover of length $m'$ of a set with chromatic number more than $c'$, then by
\ref{cresttohole}, $G$ has an odd hole of length at least $\ell$, which is impossible. Thus $G$ does not 
contain such a multicover.
Since $G$ is $(2,\phi)$-controlled by \ref{1-control}, it follows from \ref{findchand5} that $\chi(G)\le c$.
This proves \ref{mainthmagain}.~\bbox

\end{document}